\documentclass[12pt,reqno]{amsart}
\usepackage{amsmath,bm}
\usepackage{amsfonts}
\usepackage{amssymb}
\usepackage{amsthm}
\usepackage{amscd}
\usepackage{a4wide}
\usepackage{enumerate}

\usepackage{amsmath}
\usepackage{amssymb}
\usepackage{amsthm}
\usepackage{graphics}
\usepackage{eucal}
\usepackage{mathrsfs}

\usepackage{color}
\usepackage[pdftex,bookmarks,colorlinks,breaklinks]{hyperref}  

\hypersetup{linkcolor=red,citecolor=blue,filecolor=dullmagenta,urlcolor=darkblue} 
\newtheorem{theorem}{Theorem}[section]

\newtheorem{lemma}{Lemma}[section]

\newtheorem{conjecture}{Conjecture}[section]
\numberwithin{equation}{section}

\newcommand{\be}{\mathbf{e}}
\newcommand{\eo}{\mathbf{e}_{0}}

\newcommand{\en}{\mathbf{e}_{n}}

\newcommand{\eej}{\mathbf{e}_{J}}

\newcommand{\SL}{\operatorname{SL}}
\newcommand{\diag}{\operatorname{diag}}

\newcommand{\Span}{\operatorname{span}}

\newcommand{\Mat}{\operatorname{Mat}}

\newcommand{\bbR}{\mathbb{R}}

\newcommand{\bbZ}{\mathbb{Z}}
\newcommand{\bbN}{\mathbb{N}}

\newcommand{\bs}{\mathbf{s}}
\newcommand{\bq}{\mathbf{q}}

\newcommand{\bx}{\mathbf{x}}

\newcommand{\scH}{\mathcal{H}}

\newcommand{\cL}{\mathcal{L}}
\newcommand{\bp}{\mathbf{p}}

\newcommand{\bv}{\mathbf{v}}
\newcommand{\bw}{\mathbf{w}}
\newcommand{\f}{\mathbf{f}}

\newcommand{\ba}{\mathbf{a}}
\newcommand{\Id}{\operatorname{Id}}

\newcommand{\Leb}{\operatorname{Leb}}
\newcommand{\Vb}{V_{s+1 \to n}}
\newcommand{\Vd}{V_{0 \to n}}
\newcommand{\Vo}{V_{1 \to n}}

\newcommand{\bc}{\mathbf{c}}
\newcommand{\pb}{\pi_{\bullet}}
\newcommand{\oj}{\omega_{j}}

\newcommand{\scM}{\mathcal{M}}

\newcommand{\cW}{\mathcal{W}}
\newcommand{\cM}{\mathcal{M}}
\newcommand{\ff}{\mathfrak{f}}

\begin{document}
\title[Diophantine approximation on subspaces]{Diophantine approximation on subspaces of $\bbR^n$ and dynamics on homogeneous spaces}
\begin{abstract}
In recent years, the ergodic theory of group actions on homogeneous spaces has played a significant role in the metric theory of Diophantine approximation. We survey some recent developments with special emphasis on Diophantine properties of affine subspaces and their submanifolds.
\end{abstract}


\author{Anish Ghosh}
\thanks{The author acknowledges support of an ISF UGC grant.}
\address{School of Mathematics
Tata Institute of Fundamental Research
Homi Bhabha Road, Mumbai, India 400005}
\email{ghosh@math.tifr.res.in}

\maketitle
\tableofcontents

\section{Introduction}

The subject of ``metric Diophantine approximation on manifolds" is concerned with inheritance of Diophantine properties, generic with respect to Lebesgue measure in $\bbR^n$, by proper subsets.  The following ``meta conjecture"  from a survey of Kleinbock \cite{Kleinbock-survey} aptly encapsulates the general theme:
\begin{conjecture}
ÒAnyÓ Diophantine property of vectors in an ambient space (e.g. $\bbR^{n}$) which holds for almost all points in this space should hold for generic points on a nondegenerate smooth submanifold $\mathcal{M}$ of the space.
\end{conjecture}

\noindent   Let $U$ be an open subset of $\bbR^d$ and let $\f : U \to \bbR^n$ be a smooth map. Then $\f$ is said to be \emph{nondegenerate} at $\bx \in U$ if $\bbR^n$ is spanned by the
partial derivatives of $\f$ at $\bx$ of order up to $l$ for some $l$. We will call $\f$ nondegenerate if it is nondegenerate at almost every point of $U$, and the smooth manifold $\mathcal{M} = \{\f(\bx)~:~\bx \in U\}$ nondegenerate, if $\f$ is. Informally, a nondegenerate manifold is a ``sufficiently curved" manifold. An example of a nondegenerate manifold, and historically the earliest studied case is the Veronese curve:
\begin{equation}\label{ver}
(x, x^2, \dots, x^n).
\end{equation}

\noindent An example of a Diophantine property alluded to above is the following. We consider $\bx \in \bbR^n$ as a row vector, i.e. a $(1 \times n)$ matrix with real entries and $\bq \in \bbZ^n$ as a column vector. Then for every $\bx \in \bbR^n$, there exists infinitely many $\bq \in \bbZ^n$ and $p \in \bbZ$ such that
\begin{equation}\label{Dirichlet}
|p +  \bx \cdot \bq| < \|\bq\|^{-n}.
\end{equation}

\noindent Here and henceforth, $\|~\|$ refers to the supremum norm. This result is a corollary of Dirichlet's famous result which in a sense, is the beginning of Diophantine approximation. It is known that the exponent $n$ above is optimal. More precisely, say that a vector $\bx \in \bbR^n$ is \emph{very well approximable} if for some $\epsilon > 0$, there are infinitely many $\bq \in \bbZ^n$ such that
\begin{equation}\label{vwa}
|p +  \bx \cdot \bq| < \|\bq\|^{-n(1+\epsilon)}.
\end{equation}
It is a straightforward consequence of the Borel-Cantelli lemma that almost every vector in $\bbR^n$ is \emph{not} very well approximable.\\ 

Mahler \cite{Mahler} conjectured that for almost every $x \in \bbR, (x, x^2, \dots, x^n)$ is not very well approximable. In other words, the property of $\bbR^n$ of generic vectors not being very well approximable, is inherited by the curve (\ref{ver}). The difficulty of course is that a proper measurable subset of $\bbR^n$ such as the Veronese curve has zero $n$-dimensional Lebesgue measure so a priori, it could possess very different Diophantine properties. This conjecture was proved by Sprind\v{z}huk \cite{Sp1, Sp2} who in turn conjectured \cite{Sp4} that almost every vector on an analytic nondegenerate manifold is not very well approximable. Such manifolds (or the maps parametrising them) are referred to as \emph{extremal}. Sprind\v{z}huk's conjecture was proved in a stronger form by D. Kleinbock and G. Margulis in their important paper \cite{KM}. They relaxed the analyticity condition and also proved a more general \emph{multiplicative} conjecture. See \S \ref{multi} for the setup of multiplicative Diophantine approximation. Their paper is also striking because it introduced techniques from the ergodic theory of group actions on homogeneous spaces into the study of Diophantine approximation on manifolds, namely quantitative non divergence estimates for certain flows on homogeneous spaces. These techniques have since proved to be highly influential in the theory. Subsequently, the analogues of  Sprind\v{z}huk's conjectures in the $p$-adic (more generally $S$-arithmetic) setting were proved in \cite{KT} and over local fields of positive characteristic in \cite{G-pos}.\\ 

Another basic result in classical Diophantine approximation is the Khintchine-Groshev theorem which provides a very satisfactory measure theoretic generalisation of Dirichlet's theorem. Let $\psi : \bbR_{+} \cup \{0\} \to \bbR_{+} \cup \{0\}$ be a non-increasing function and define $\cW(\bbR^n, \psi)$ to be the set of $x \in [0,1]^n \subset\bbR^n$ for which there exist infinitely many $\bq \in \bbZ^{n}$ such that
\begin{equation}\label{preKG}
|p +  \bx \cdot \bq| < \psi(\|\bq\|^n)
\end{equation}
\noindent for some $p \in \bbZ$. The Khintchine-Groshev Theorem (\cite{Khintchine}, \cite{Groshev}, and see \cite{Dodson} for a nice survey) gives a characterization of the measure of $\cW([0, 1]^n, \psi)$ in terms of $\psi$:

\begin{theorem}\label{KG}
\begin{equation} 
\Leb(\cW([0, 1]^n, \psi)) = \left\{ 
\begin{array}{rl} 
0 & \text{if } \sum \psi(k) < \infty\\ 
\\
1 & \text{if } \sum \psi(k) = \infty.
\end{array} \right.
\end{equation}
\end{theorem}

\noindent Here $\Leb$ denotes Lebesgue measure. In \cite{BKM}, Bernik, Kleinbock and Margulis  proved that the convergence case of the Khintchine-Groshev theorem holds for nondegenerate manifolds in both standard and multiplicative cases. The standard case was independently proved by Beresnevich in \cite{Ber1}. It should be noted that the convergence case of the Theorem for vectors in $\bbR^n$ is once more a straightforward consequence of the Borel-Cantelli lemma but the proof for manifolds is substantially more complicated. The complementary divergence case of the Theorem for nondegenerate manifolds was established in \cite{BBKM}. We refer the reader to \cite{Kleinbock-diophmeas} for a nice survey of these and related developments including inheritance of other Diophantine properties by nondegenerate manifolds like the paucity of singular and Dirichlet improvable vectors.\\  

\noindent The present survey is about the opposite end of the spectrum, namely the ``non-nondegenerate" manifolds. An example is an affine subspace and it is with Diophantine approximation on affine subspaces and their nondegenerate submanifolds that we shall be primarily concerned. The concept of nondegeneracy in an affine subspace is a natural extension of non degeneracy in $\bbR^n$ and a definition is provided below.  The meta conjecture above excludes affine subspaces for a good reason because it is possible to construct examples, say lines in $\bbR^n$ whose slopes are rational or very close to rational and which do not possess this inheritance phenomena for generic Diophantine properties. For example, not all lines are extremal. On the other hand, it can be shown that a large class of affine subspaces do possess inheritance properties. The study of Diophantine approximation on affine subspaces therefore provides a highly interesting interplay between Diophantine properties of the matrix parametrising the subspace and the generic Diophantine properties it inherits. 

\subsection*{Acknowledgement} The author thanks the referee for helpful comments which have helped improve the exposition.

\section{Diophantine approximation on affine subspaces and their submanifolds}

\noindent The main questions as regards Diophantine approximation on affine subspaces are as follows:

\begin{enumerate}
\item Under which conditions are Diophantine properties which are generic for Lebesgue measure inherited by affine subspaces and their nondegenerate submanifolds?\\
\item Is it possible to explicitly describe the subspaces with the inheritance property in terms of Diophantine properties of their parametrising matrices?
\end{enumerate}

\noindent To properly address these issues we need to describe the Diophantine properties of parametrising matrices of affine subspaces which come into play. For $v > 0$ and $m, n \in \bbN$, denote by $\cW_{v}(m,n)$ the set of matrices 
$A \in \Mat_{m, n}$ for which there are infinitely many $\bq \in \bbZ^n$ such that
\begin{equation}\label{defexponent}
\|A \bq + \bp\| < \|\bq\|^{-v}
\end{equation}
\noindent for some $\bp \in \bbZ^m$. Then $ \cW_{v_2}(m, n) \subset \cW_{v_2}(m, n) \text{ if } v_1 \leq v_2$. Following \cite{Kleinbock-extremal}, we define
\begin{equation}
 \cW^{+}_{v}(m, n) := \bigcup_{u > v}  \cW_{u}(m, n) ~\text{ and }~  \cW^{-}_{v}(m, n) := \bigcap_{u < v}  \cW_{u}(m, n).
\end{equation}

\noindent By Dirichlet's theorem for matrices, $ \cW_{n/m}(m, n) = \Mat_{m, n}(\bbR)$ and $ \cW^{+}_{n/m}(m, n)$ is the set of very well approximable matrices.\\

The Diophantine exponent $\omega(A)$ of a matrix $A \in \Mat_{m \times n}(\bbR)$ is defined to be the supremum of $v > 0$ for which there are infinitely many $\bq \in \bbZ^n$ such that
\begin{equation}\label{defexponent}
\|A \bq + \bp\| < \|\bq\|^{-v}
\end{equation}
\noindent for some $\bp \in \bbZ^m$. It is well known that $n/m \leq \omega(A) \leq \infty$ for all $A \in \Mat_{m \times n}(\bbR)$ and that $\omega(A) = n/m$ for Lebesgue almost every $A$. The exponent is closely related to the sets $\cW_{v}(m,n)$ of course, 
$$ A \notin \cW^{+}_{n}(n, 1) \text{ if and only if } \omega(A) \leq n.
 $$

\noindent The first known result in the subject is due to W. M. Schmidt \cite{Schmidt} who studied extremality for lines in $\bbR^n$.

\begin{theorem}\label{Schmidt}\cite{Schmidt}
Let $\mathcal{L}$ be a line in $\bbR^n, n \geq 2, \mathcal{L} = (t, \alpha_1t+\beta_11, \dots , \alpha_{n-1}t+ \beta_{n-1})$, where either 
\begin{equation}
\omega(\alpha_1, \dots, \alpha_{n-1}) = n - 1 \text{ or } \omega(\beta_1,\dots, \beta_{n-1}) = n-1.
\end{equation}
Then $\mathcal{L}$ is extremal.
\end{theorem}

\noindent Subsequently Beresnevich, Bernik, Dickinson and Dodson \cite{BBDD} studied the Khintchine-Groshev theorem for lines in $\bbR^n$. Let us say that a smooth manifold is Groshev type for convergence (resp. Groshev type for divergence) if the convergence (resp. divergence) case of the Khintchine-Groshev theorem holds for it.

\begin{theorem}\label{BBDD}(Beresnevich-Bernik-Dickinson-Dodson)
Let $\mathcal{L}$ be a line passing through the origin, parametrised by a matrix $A$. Then $\mathcal{L}$ is Groshev type for convergence, as well as Groshev type for divergence whenever $\omega(A) < n$.
\end{theorem}

A systematic study of Diophantine approximation on subspaces and their submanifolds was initiated by D. Kleinbock in \cite{Kleinbock-extremal}. To describe the results in this paper we need the notion of nondegeneracy of manifolds in affine subspaces. Let $\scH$ be an affine subspace of $\bbR^n$, let $U$ be an open subset of $\bbR^d$ and let $\f : U \to \bbR^n$ be a differentiable map. Then $\f$ is said to be nondegenerate in an affine subspace $\scH$ of $\bbR^n$ at $x_0 \in U$ if $\f(U) \subset \scH$ and the span of all the partial derivatives of $\f$ up to some order is the linear part of $\scH$. Let $\mathcal{M}$ be a $d$-dimensional submanifold of $\mathcal{H}$. We call $\mathcal{M}$ nondegenerate in $\scH$ at $y \in \scM$ if any (equivalently, some) diffeomorphism $\f$ between an open subset $U$ of $\bbR^d$ and a neighborhood of $y$ in $\scM$ is nondegenerate in $\scH$ at $\f^{-1}(y)$. We say that $\f : U \to \scH$ is nondegenerate in $\scH$ if it is nondegenerate at almost every point of $U$. Equivalently we say that $\mathcal{M}$ is nondegenerate in $\mathcal{L}$ if it is nondegenerate in $\mathcal{L}$ at almost every point $y \in \mathcal{M}$ with respect to the natural measure class on $\mathcal{H}$.  In \cite{Kleinbock-extremal}, Kleinbock proved the following two beautiful theorems resolving the analogue of Sprind\v{z}uk's conjecture for affine subspaces and their nondegenerate submanifolds.

\begin{theorem}\label{Kleinbock1}(\cite{Kleinbock-extremal} Theorem $1.2$)
Let $\scH$ be an affine subspace of $\bbR^n$. Then:
\begin{enumerate}[(a)]
\item if $\scH$ is extremal and $\f : U \to \scH, U \subset \bbR^d$, is a smooth map which is nondegenerate in $\scH$, then $\f$ is extremal; 
\item if $\scH$ is not extremal, then all points of $\scH$ are very well approximable (in particular, no subset of $\scH$ is extremal).
\end{enumerate}

\end{theorem}

\noindent Part $(ii)$ in the above theorem presents a striking dichotomy as regards very well approximable vectors on affine subspaces. It would be of considerable interest to investigate if such a dichotomy exists for other Diophantine properties. 

In fact, sufficient and sometimes necessary conditions for affine subspaces to be extremal in terms of the matrix parametrising the subspace can also be found. 

Let $\scH$ be an $s$ dimensional affine subspace of $\bbR^n$. Without loss of generality, one can choose a parametrizing map of the form $\bx \to (\bx, \bx A' + a_0)$, where $A'$ is a matrix of size $s \times (n-s)$ and $a_0 \in \bbR^{n-s}$. Denote the vector $(1, x_1, \dots, x_s)$ by $\tilde{\bx}$, and
the matrix $\begin{pmatrix}\ba_0\\A' \end{pmatrix}  $ by $A \in \Mat_{s+1, n-s}$. Then $\scH$ is parametrised by the map 
\begin{equation}\label{defaff}
\bx \to (\bx, \tilde{\bx}A). 
\end{equation}
In \cite{Kleinbock-extremal}, D. Kleinbock provided a necessary and sufficient condition on $A$ for extremality  of the subspace $\scH$ parametrised as above in the following two cases:
\begin{enumerate}
\item $\scH$ is an affine hyperplane, i.e. has dimension $s = n - 1$;\\
\item $\scH$ is a line passing through the origin in $\bbR^n$.
\end{enumerate}

\noindent He proved

\begin{theorem}\label{Kleinbock2}(\cite{Kleinbock-extremal} Theorem $1.3$)
In the two cases above, the map (\ref{defaff}) is extremal if and only if 
\begin{equation}\label{diocond1}
A \notin \cW^{+}_{n}(s+1, n - s).
\end{equation}
\end{theorem}

\noindent The method of proof for the two cases is different. For the hyperplane case, a modification of the dynamical approach of Kleinbock-Margulis is used while for the line passing through the origin, a modification of the method of \cite{BBDD}. In \cite{G1, G-div} the Khintchine-Groshev theorem was investigated for affine hyperplanes.
\begin{theorem}\label{G1}(\cite{G1, G-div})
Let $\scH$ be an affine hyperplane parametrised as in (\ref{defaff}). Assume that
\begin{equation}\label{std}
A \notin \mathcal{W}_{n}^{-}(n, 1).
\end{equation}
\noindent Then $\scH$ is Groshev type for convergence (resp. divergence), namely, the set
\begin{equation}
\{\bx \in \bbR~:~| (\bx, \tilde{\bx}A) \cdot \bq + p| < \psi(\|\bq\|) \text{ for infinitely many } \bq \in \bbZ^n, p \in \bbZ\}
\end{equation}
\noindent  has zero (resp. full) Lebesgue measure whenever 
\begin{equation}
\sum_{k=1}^{\infty}k^{n-1}\psi(k)
\end{equation}
\noindent converges (resp. diverges).
\end{theorem}

\noindent The set of matrices satisfying (\ref{std}) is large. Indeed, using Dodson's formula \cite{Dodson} for Hausdorff dimension, one can show that the set of matrices not satisfying (\ref{std}) has Hausdorff dimension $1$. However, obtaining a necessary condition for an affine subspace to be of Groshev type for convergence and or divergence seems to be a difficult open problem. We note that neither the result in \cite{BBDD} nor the Theorem above provides a necessary condition. It is plausible that given a fixed function $\psi$ that such a condition may be possible to compute. The method of proof in \cite{G1} uses a modification of a technique due to Bernik, Kleinbock and Margulis and uses non divergence estimates for certain flows on homogeneous spaces, while the proof of the divergence counterpart uses estimates from \cite{G1} in addition to the method of ``regular systems".

\section{Higher Diophantine exponents of matrices}
An upshot of the previous section is that in certain cases, e.g. for certain lines and for co-dimension one subspaces, it is relatively easy to describe the properties of affine subspaces needed in order for them to inherit extremality or Groshev type behaviour. In order to explore Diophantine properties of higher co-dimension affine subspaces and their nondegenerate submanifolds, we need the notion of higher Diophantine exponent. Let $\scH$ be an $s$ dimensional affine subspace of $\bbR^n$. We take $U \subset \bbR^d$ and $\f = (f_1, f_2, \dots, f_n) : U \to \scH$ a smooth nondegenerate map. As in the previous section, $\scH$ can be parametrized by the map (\ref{defaff}). By our assumptions above, we have that $1, f_1, \dots, f_s$ are linearly independent over $U$. Setting $\ff = (f_1, f_2, \dots, f_s)$, and $\tilde{\ff} = (1, f_1, \dots, f_s)$ we can write $\f$ as 
\begin{equation}\label{frep}
x \to (\ff(x), \tilde{\ff}(x)A)
\end{equation}

Following Kleinbock \cite{Kleinbock-exponent}, we now define the higher Diophantine exponents of $A$. Though their definition is not as transparent, as shown in \cite{Kleinbock-exponent}, they encode important Diophantine properties of the affine subspace parametrized by $A$. In order to motivate the definition of these, we begin, as in \cite{G-Monat}, by considering Khintchine-type inequalities on affine sub spaces and their submanifolds. In view of the notation above, this amounts to studying  the inequality
\begin{equation}\label{reform}
|p + (\ff(x), \tilde{\ff}(x)A)\bq| < \psi(\|\bq\|^n)
\end{equation}
\noindent For $A \in \Mat_{s+1 \times n-s}(\bbR)$, we set 
\begin{equation}\label{defR}
R_A = \begin{pmatrix}\Id_{s+1} & A \end{pmatrix}.
\end{equation}
\noindent Then (\ref{reform}) can be written as
\begin{equation}
|\tilde{\ff}(x)R_{A}\begin{pmatrix}p \\ \bq\end{pmatrix}| < \psi(\|\bq\|^n).
\end{equation}
\noindent The matrix $R_A$ plays an important role in the Diophantine properties of $\f$. This role is expressed in the form of certain ``higher" Diophantine exponents. Let $\eo, \dots, \en$ denote the standard basis of $\bbR^{n+1}$ and set 
\begin{equation}
V_{i \to j} = \Span(\be_i, \dots, \be_j).
\end{equation}
\noindent Let $\bw \in \bigwedge^{j}(\Vd)$ represent a discrete subgroup $\Gamma$ of $\bbZ^{n+1}$. Define the map $\bc : \bigwedge^{j}(\Vd) \to (\bigwedge^{j-1}(\Vo))^{n+1}$ by
\begin{equation}\label{defc}
\bc(\bw)_i = \sum_{\substack{J \subset \{1, \dots, n\}\\ \#J = j-1}} \langle \be_i \wedge \eej, \bw\rangle \eej
\end{equation} 
\noindent and let $\pb$ denote the projection $\bigwedge(\Vd) \to \bigwedge(\Vb)$. For each $j = 1, \dots, n-s$, define
\begin{equation}\label{defexponenthigher}
\oj(A) = \sup\left\{ v\left| \aligned \exists\, 
\bw \in  \bigwedge^{j}(\bbZ^{n+1}) \text{ with arbitrary large
   } \|\pb(\bw)\|  \\
\text{   such that   }\|R_{A}
\bc(\bw)\|
< \|\pb(\bw)\|^{-\frac {v+1-{j}}{j}}\ \ 
\endaligned\right.\right\}
\end{equation}
\noindent It follows from Lemma 5.3 in loc. cit. that $\omega_1(A) = \omega(A)$ thereby justifying the terminology. Computing higher Diophantine exponents is rather difficult, we refer the reader to \cite{Kleinbock-exponent} and \cite{G-Monat}. In \cite{Kleinbock-exponent}, D. Kleinbock proved that the extremality of an affine subspace is characterised by it's higher Diophantine exponents.

\begin{theorem}[Corollary $5.2$ and Theorem $0.3$ in \cite{Kleinbock-exponent} and Theorem $1.2$ in \cite{Kleinbock-extremal}]\label{kleinbock}
Let $\scH$ be an $s$-dimensional affine subspace parametrized as in (\ref{defaff}). 
\begin{enumerate}
\item $\scH$ is extremal if and only if 
\begin{equation}\label{Kl-condition}
\omega_{j}(A) \leq n ~\text{for every}~j = 1, \dots, n-s.
\end{equation}
\item If $\scH$ is extremal, then so is  $\scM$ for any smooth nondegenerate submanifold $\scM$ of $\scH$.
\end{enumerate} 
\end{theorem}
 
\noindent Using the transference principle of Beresnevich and Velani \cite{BV}, the Theorem above implies an \emph{inhomogeneous} version of extremality for affine subspaces and their submanifolds. We refer the reader to \cite{BV} for details on the inhomogeneous setup.  In \cite{G-Monat}, we proved that weakening the condition (\ref{Kl-condition}) suffices to prove the convergence case of the Khintchine-Groshev theorem, thereby confirming a conjecture of Kleinbock, cf. \cite{Kleinbock-exponent} \S $6.5$.

\begin{theorem}\label{main}
Let $\scH$ be an $s$-dimensional affine subspace parametrized as in (\ref{defaff}). Assume that 
\begin{equation}\label{diocond}
\oj(A) < n~\text{for every}~j = 1, \dots, n-s. 
\end{equation} 
\noindent Then
\begin{enumerate}
\item $\scH$ is Groshev type for convergence.\\
\item $\scM$ is Groshev type for convergence for every smooth nondegenerate submanifold $\scM$ of $\scH$.
\end{enumerate}
\end{theorem}

\noindent The proof follows the dynamical approach once more. While the corresponding divergence case has not been proved yet, we expect that the use of regular systems as in \cite{BBKM, G-div} should yield the divergence case of the Khintchine-Groshev theorem for affine subspaces and their submanifolds under the same conditions, namely (\ref{diocond}). As in the case of hyperplanes, a necessary condition seems elusive at present.\\

\section{Variations: Multiplicative and weighted Diophantine approximation}\label{multi}
 For an integer $n$ and a vector $\bq = (q_1, \dots, q_n) \in \bbR^{n}$, set
\begin{equation}\label{normsdef}
\Pi_{+}(\bq) := \prod_{i = 1}^{n}\max(|q_i|, 1).
\end{equation}
\noindent The study of Diophantine approximation with the supremum norm replaced by \ref{normsdef} is referred to as \emph{Multiplicative} Diophantine approximation. This replacement leads to many subtle difficulties in metric theory. We refer the reader to \cite{Bugeaud} for a survey of some problems in the subject. 

\noindent One can then prove a multiplicative version of Dirichlet's theorem and of the Khintchine-Groshev theorem and define the set of \emph{very well multiplicatively approximable} vectors in $\bbR^n$ which also form a set of measure zero. In \cite{KM}, the stronger (than extremality) property that almost every point on a nondegenerate manifold is not very well multiplicatively approximable \footnote{such a manifold is referred to as \emph{strongly extremal}} was established, thereby resolving a conjecture of Baker  and in \cite{BKM}, the multiplicative  of the Khintchine-Groshev theorem for nondegenerate manifolds were established. In the setting of affine subspaces, Kleinbock proved multiplicative analogues of all his results for standard Diophantine approximation. For instance,

\begin{theorem}\label{Kleinbock2}(\cite{Kleinbock-extremal} Theorem $1.4$)
Let $\scH$ be an affine subspace of $\bbR^n$. Then:
\begin{enumerate}[(a)]
\item if $\scH$ is strongly extremal and $\f : U \to \scH, U \subset \bbR^d$, is a smooth map which is nondegenerate in $\scH$, then $\f$ is strongly extremal; 
\item if $\scH$ is not strongly extremal, then all points of $\scH$ are very well multiplicatively approximable.
\end{enumerate}

\end{theorem}

\noindent Similarly, in \cite{G-mult}, a multiplicative form of Khintchine's theorem for affine hyperplanes is proved. For a column matrix $A = (a_i)$, set

\begin{equation}\label{notzero}
r = r(A) = \#\{1 \leq i \leq n-1~|~a_i \neq 0\}.
\end{equation}

\noindent Then

\begin{theorem}\label{main}
Assume that $r(A) = n-1$ and that $A$ satisfies (\ref{std}). Then
\begin{equation}
\{\bx \in B~:~| (\bx, \tilde{\bx}A) \cdot \bq + p| < \psi(\Pi_{+}(\bq)) \text{ for infinitely many } \bq \in \bbZ^n, p \in \bbZ\}
\end{equation}
\noindent has zero measure whenever
$$\sum_{k = 1}^{\infty}(\log k)^{n-1}\psi(k) < \infty.$$
\end{theorem}

Another possible modification is to consider a \emph{weighted} norm. For $\bs \in \bbR^{n}_{+}$ such that $\sum_{i = 1}^{n}s_i = 1$, we follow \cite{BKM} in defining the $\bs$-quasinorm on $\bbR^n$ by
\begin{equation}\label{defsnorm}
\|\bx\|_{\bs} := \max_{1 \leq i \leq n}|x_i|^{1/s_i}.
\end{equation}

\noindent One  could then consider the corresponding Diophantine inequalities. A strength of the dynamical approach to Diophantine approximation outlined in the last section is that  it is flexible enough to address all these variations.

\section{Diophantine exponents of measures}
\noindent Following \cite{Kleinbock-exponent}, for a Borel measure $\mu$ we define its Diophantine exponent by 
\begin{equation}
\omega(\mu) = \sup\{v ~|~ \mu\{\bx~|~ \omega(\bx) > v\} > 0\}.
\end{equation}

\noindent The definition only depends on the measure class of $\mu$. Recall that $\lambda$ denotes Lebesgue measure on $\bbR^n$. If $\mathcal{M}$ is a smooth submanifold of $\mathbb{R}^n$ parametrised by a smooth map $\f$, then set the Diophantine exponent $\omega(\mathcal{M})$ to be equal to $\omega(\f_{*}\lambda)$ where $\f_{*}\lambda$ is the push forward of $\lambda$ by $\f$. Then $\omega(\mu)$ is always at least equal to $n$ and $\omega(\lambda) = n$. Extremal manifolds are precisely those for which $\omega(\f_{*}\lambda) = n$. In \cite{Kleinbock-exponent}, D. Kleinbock showed that Diophantine exponents of affine subspaces are inherited by their nondegenerate submanifolds, thereby generalising his results on extremality.

\begin{theorem}( \cite{Kleinbock-exponent} Theorem $0.3$)
Let $\mathcal{L}$ be an affine subspace of $\mathbb{R}^n$, and let $\mathcal{M}$ be a submanifold of $\mathcal{L}$ which is nondegenerate in $\mathcal{L}$. Then
$$\omega(\mathcal{M}) = \omega(\mathcal{L}) = \inf\{\omega(\bx)~|~ \bx \in \mathcal{L}\} = \inf\{\omega(\bx)~|~ \bx \in \mathcal{M}\}.$$

\end{theorem}

\noindent Further, in the same paper, Kleinbock computes the Diophantine exponent of a subspace in terms of the exponent of its parametrising matrix. 
\begin{theorem}\label{Kleinbock-exp} Let $\mathcal{L}$ be an affine subspace of $\bbR^n$ parametrised by a matrix $A$ such that either (a) all the columns or (b) all the rows of $A$ are rational multiples of one column (resp. one row). Then,
\begin{equation} 
\omega(\mathcal{L}) = \max(\omega(A), n).
\end{equation}
\end{theorem}
\noindent Similarly, the multiplicative exponent of a vector can be defined as:
\begin{equation}\label{multexp}
\omega^{\times}(\bx) := \sup\left\{v~|~\exists \text{ infinitely many } \bq \in \bbZ^n \text{ such that } |\bx \cdot \bq + p| \leq \Pi_{+}(\bq)^{-v/n} \text{ for some } p \in \bbZ \right\}
\end{equation}

\noindent and that of a measure can be defined as

\begin{equation}\label{multexp}
\omega^{\times}(\mu) := \sup\{v~|~\mu\{\bx~|~\omega^{\times}(\bx) > v\} > 0\}.
\end{equation}

\noindent In \cite{Zhang2}, Y. Zhang calculates the multiplicative Diophantine exponent of affine hyperplanes and their nondegenerate submanifolds. Here is Theorem $1.4$ from \cite{Zhang2}:

\begin{theorem}\label{zhang2}
If $\cL$ is a hyperplane of $\bbR^{n}$ and $\cM$ is a nondegenerate submanifold in $\cL$, then
$$\omega^{\times}(\cL) = \omega^{\times}(\cM) = \inf \{\omega^{\times}(\bx)~|~ \bx \in \cL\} = \inf\{\omega^{\times}(\bx)~|~\bx \in \cM\}. $$
\end{theorem}

\noindent She further explicitly computes the exponents in terms of the parametrising coefficients of the affine hyperplane. In \cite{Zhang2} it is proved that

\begin{theorem}\label{zhang3}
Suppose the affine hyperplane $\cL$ is parametrised by the matrix $A$ and let $r=r(A)$ be as in (\ref{notzero}). Then
 $$\omega^{\times}(\cL) = \max\left(n, \frac{n}{r+1}\sigma(A)\right).$$
\end{theorem}

\noindent Here
$$\sigma(A):= \sup\{v~|~\exists \text{ infinitely many } q \in \bbZ \text{ with } \|qA+ \bp\|< |q|^{-v} \text{ for some } \bp \in \bbZ \}$$

\noindent is the simultaneous Diophantine exponent of $A$.  The  simultaneous and linear Diophantine exponent encountered before are related to each other via Khintchine's transference principle, we refer the reader to \cite{Khintchine-transfer, Cassels}.

\section{Badly approximable vectors}
A real number $x$ is called \emph{Badly approximable} if there exists a constant $c = c(x) > 0$ such that
\begin{equation}\label{def:bad1}
\left| x - \frac{p}{q}\right| \geq \frac{c}{q^2}
\end{equation}
\noindent for all $(p, q) \in \bbZ \times \bbN$. Badly approximable numbers are precisely those whose whose partial quotients are bounded. They form a set of zero Lebesgue measure but of full Hausdorff dimension. In higher dimensions, one can adopt either the simultaneous point of view or the ``dual" point of view. We adopt the dual perspective. Then a vector $\bx \in \bbR^n$ is badly approximable if there exists a constant $c(\bx) > 0$ such that
\begin{equation}\label{def:bad2}
|\bx\cdot\bq + p| \geq \frac{c}{\|\bq\|^n}
\end{equation}

\noindent for all $(\bq, p) \in \bbZ^n \times \bbZ$. Let $\mathbf{Bad}$ be the set of badly approximable vectors. Then $\mathbf{Bad}$ has full Hausdorff dimension and zero Lebesgue measure. In recent years there has been significant progress in understanding the structure of badly approximable vectors. In an important work, V. Beresnevich \cite{Ber2} resolved a long standing problem posed by H. Davenport \cite{Davenport} on the size of badly approximable vectors on manifolds.  

\begin{theorem}\label{Beresnevich}(Beresnevich)
Let $\f = (f_1,\dots, f_n) : B \to \bbR^n$ be an analytic map defined on a ball $B \subset \bbR^m$ such that the functions 
$$1, f_1, \dots, f_n \text{ are linearly independent over } \bbR.$$ 
Let $\mathcal{M}$ be the analytic manifold parametrised by $\f$. Then  $\dim (\mathcal{M} \cap \mathbf{Bad}) = \dim \mathcal{M}$. 
\end{theorem}

\noindent We note that the linear independence condition above is equivalent to non degeneracy for analytic functions. In fact, Beresnevich proves much more. He considers the set of weighted badly approximable numbers and permits intersections of finite collections of weighted badly approximable vectors with analytic nondegenerate manifolds thereby resolving in a strong form, a generalisation of W. M. Schmidt's conjecture on intersections of sets of weighted badly approximable vectors which was resolved by Badziahin, Pollington and Velani. We refer the reader to \cite{Ber2} for further details. In \cite{BGV}, we investigate badly approximable vectors on affine subspaces and their submanifolds. As in \cite{Ber2}, we work in the weighted setting but restrict ourselves to a special case in this survey for ease of exposition. 

\begin{theorem}\label{main-Bad}
Let $\scH$ be an $s$-dimensional affine subspace parametrized as in (\ref{defaff}). Assume that 
\begin{equation}\label{diocond}
\oj(A) < n~\text{for every}~j = 1, \dots, n-s. 
\end{equation} 
\noindent Then
\begin{enumerate}
\item $\dim (\mathbf{Bad} \cap \scH) = \dim \scH = s$.\\
\item $\dim (\mathbf{Bad} \cap \scM) = \dim \scM$ for every analytic nondegenerate submanifold $\scM$ of $\scH$.
\end{enumerate}
\end{theorem}

\section{Diophantine approximation with restricted numerator and denominator}

In this section, we discuss some recent progress on a more refined problem of Diophantine approximation with restrictions on the rationals. Such problems have a long history, for instance, the problem of approximating irrational numbers with rationals whose denominators are prime was studied by Vinogradov, Heath-Brown and Jia \cite{HeJi} and others. The current best possible result is due to Matom\"{a}ki \cite{Mato} and states that for irrational $\alpha$, there exist infinitely many primes $p$ such that
\begin{equation}\label{eq:res}
|\langle \alpha p \rangle| < p^{-\nu}
\end{equation}
\noindent for any $\nu < 1/3$. In this section, $\langle~\rangle$ denotes the distance to the nearest integer. The corresponding problem with both numerator and denominator prime was studied by Ramachandra and others \cite{Rama, Srini}. More generally, one can study analogues of Khintchine's theorem with restrictions. The case of prime denominator follows from the theorem of Duffin-Schaeffer \cite{DuSch} whereas Harman \cite{HarmK} has established the complete analogue of Khintchine's theorem in the case of both prime numerator and denominator. His work was extended to simultaneous approximation in higher dimensions by Jones \cite{Jones}. 
 
In \cite{HaJo}, Harman and Jones  initiated the study of Diophantine approximation on manifolds with primality constraints and proved 

\begin{theorem}(Harman-Jones) Let $\varepsilon>0$ and $\tau>1$. Then for almost all positive $\alpha$, with respect to the Lebesgue measure, there are infinitely many $p,q,r$ prime such that
\begin{equation}
0 < p\alpha - r< p^{-3/16+\varepsilon} \quad \mbox{and} \quad 0<p\alpha^{\tau} - q< p^{-3/16+\varepsilon}.
\end{equation}
\end{theorem}

Since $\tau > 1$ the curve $(\alpha, \alpha^{\tau})$ is a nondegenerate curve. The generic best possible exponent above is $-1/2$. In \cite{BG1, BG2} the study of Diophantine approximation on affine subspaces was initiated using the technique of exponential sums. Assume that ${\bf c} \in \bbR^d$ satisfies \footnote{ such a vector is usually called $k$-Diophantine}: there exists a constant $C>0$ such that
\begin{equation} \label{condi}
|\langle{\bf v}\cdot {\bf c}\rangle|>\frac{C}{||{\bf v}||^k} \quad \mbox{for every } {\bf v}\in \mathbb{Z}^d\setminus \{{\bf 0}\}.
\end{equation}

The main result in \cite{BG2} is

\begin{theorem} \label{super}
Let $d$ be a positive integer and $k\ge d$ be a positive real number. Define
\begin{equation} \label{gammadef}
\gamma_{d,k}:=\frac{1}{d(3k+2)}
\end{equation}
and suppose that $0<\varepsilon<\gamma_{d,k}$.
Let $c_1,...,c_d$ be positive irrational numbers such that the vector ${\bf c}=(c_1,...,c_d)$ satisfies (\ref{condi}). Then for almost all positive real $\alpha$, with respect to the Lebesgue measure, there are infinitely many $(d+2)$-tuples $(p,q_1,\dots,q_d,r)$ with $p$ and $r$ prime and $q_1,\dots,q_d$ positive integers such that simultaneously
\begin{equation} \label{simultan}
\begin{split}
0 & <p\alpha-r<p^{-\gamma_{d,k}+\varepsilon},\\
0 & <pc_i\alpha-q_i<p^{-\gamma_{d,k}+\varepsilon} \mbox{ for all } i\in \{1,\dots,d\}.
\end{split}
\end{equation}
\end{theorem}

\noindent In \cite{BG1}, the special case of Diophantine approximation on lines in $\bbR^2$ with primality constraints was considered earlier. It is an interesting but challenging problem to generalise Theorem \ref{super} to higher dimension subspaces and their nondegenerate submanifolds. Improving the exponent, even conditionally, is also a challenging problem.

\section{Dynamics of group actions }
We very briefly discuss the approach to proving many of the results in this survey\footnote{Other than those related to prime Diophantine approximation where analytic techniques are the most successful.}, and in particular extremality and the Khintchine-Groshev Theorem for affine subspaces.  By now, there are several very nice and comprehensive surveys on the interaction between homogeneous dynamics and number theory, we refer the reader to \cite{Kleinbock-survey, Kleinbock-Clay} in particular where connections of quantitative recurrence estimates for certain orbits on homogeneous spaces with metric Diophantine approximation, are discussed in detail. Following Dani \cite{Dani} it was realised that Diophantine properties of vectors can be encoded in the dynamics on certain group actions on the homogeneous space $\SL(n, \bbR)/\SL(n, \bbZ)$. The latter is a non compact finite volume space and can be identified with the space of unimodular lattices $X_n$ in $\bbR^n$ via the standard action of $\SL(n, \bbR)$. Kleinbock and Margulis, in their influential paper \cite{KM}, introduced dynamical methods into the study of metric Diophantine approximation on manifolds. Let 
\begin{equation}\label{def:delta}
\delta(\Lambda) := \inf_{\bv \in \Lambda \backslash \{0\}} \|\bv\|
\end{equation}
denote the length of the shortest vector of a lattice. According to Mahler's compactness criterion, compact subsets of $X_n$ are precisely those sets of lattices where $\delta$ is uniformly bounded from below. Given $\bx \in \bbR^n$, associate to it the lattice:
$$\Lambda_{\bx}:= \begin{pmatrix} 1 & \bx\\ 0 & I \end{pmatrix}\bbZ^{n+1}, $$
\noindent here $I$ denotes the identity $n \times n$ matrix. For $t \in \mathbb{Z}_{+},$ let $g_t := \diag(e^{nt}, e^{-t}, \dots, e^{-t})$. Then it is shown in \cite{KM} that 
\begin{lemma}
The vector $\bx$ is very well approximable if and only if there exists $\gamma > 0$ and infinitely many $t \in \mathbb{Z}_{+}$ such that
$$ \delta(g_t\Lambda_{\bx}) \leq e^{-\gamma t}.$$
\end{lemma}
\noindent This translates the problem of studying very well approximable vectors into one of understanding orbits of group actions on the space of lattices and the main tool used to prove Sprind\v{z}uk's conjecture is a quantitative non divergence estimate for certain polynomial like flows on $X_n$. Namely, the fact that almost every vector in $\bbR^n$ is not  very well approximable, translates to the fact that for almost every lattice, certain orbits recur in a quantitative manner, to compact sets. 
We refer the reader to \cite{Kleinbock-Clay} for the history of non divergence estimates for unipotent flows and for a detailed proof of the aforementioned quantitative estimate. The proof of the convergence case of the Khintchine-Groshev theorem for nondegenerate manifolds in \cite{BKM} and the subsequent proof for affine subspaces in \cite{G1, G-Monat} and other works uses a variation on this result as the primary technical tool. The nondegeneracy of the maps (reps. the Diophantine conditions assumed on affine subspaces) arise naturally in this dynamical reformulation as obstructions to ``diverging to infinity" in the space of lattices.

\bibliographystyle{amsplain}

\end{document}